\magnification 1200
\def\R{{\rm I\kern-0.2em R\kern0.2em \kern-0.2em}}
\def\N{{\rm I\kern-0.2em N\kern0.2em \kern-0.2em}}
\def\P{{\rm I\kern-0.2em P\kern0.2em \kern-0.2em}}
\def\B{{\rm I\kern-0.2em B\kern0.2em \kern-0.2em}}
\def\Z{{\rm I\kern-0.2em Z\kern0.2em \kern-0.2em}}
\def\C{{\bf \rm C}\kern-.4em {\vrule height1.4ex width.08em depth-.04ex}\;}
\def\B{{\bf \rm B}\kern-.4em {\vrule height1.4ex width.08em depth-.04ex}\;}

\def\D{{\Delta}}
\def\DD{{\overline \Delta}}
\def\bD{{b\Delta}}
\def\z{{\zeta}}
\def\cC{{\cal C}}

\def\L{{\cal L}}
\def\cC{{\cal C}}
\def\cS{{\cal S}}
\def\cS{{\cal S}}

\def\cT{{\cal T}}
\def\F{{\cal F}}

\def\Tt{{(T-t)(T-1/t)}}
\font\ninerm=cmr8
\
\vskip 25mm
\centerline {\bf SMALL FAMILIES OF COMPLEX LINES FOR TESTING}

\centerline {\bf HOLOMORPHIC EXTENDIBILITY}
\vskip 6mm
\centerline{Josip Globevnik}
\vskip 6mm
{\noindent\ninerm \hfill To Urban, Manja and Lidija  }
\vskip 6mm
{\noindent \ninerm ABSTRACT \ \ Let $B$ be the open unit ball in $\C^2$ and let
$a, b\in\overline B,\ a\not=b$.
It is known that given $k\in\N$ there is a function
$f\in C^k(bB)$ which extends holomorphically into $B$ along
any complex line passing through either $a$ or $b$ yet $f$
does not extend holomorphically through $B$. In the paper we
show that there is no such function in $C^\infty (bB)$. Moreover, we obtain a
fairly complete description of pairs of points $a, b\in\C^2, a\not= b,$ such that if
$f\in C^\infty (bB)$ extends holomorphically into $B$ along every complex line passing
through either $a$ or $b$ that meets $B$, then $f$ extends holomorphically through $B$.}
\vskip 6mm
\bf 1.\ Introduction and the main result\rm
\vskip 2mm Denote by $\D$ the open unit disc in $\C $ and by $B$ the open unit ball in $\C^2$. If
$f$ is a continuous function on $bB$ and $L$ is a complex line in $\C^2$ that meets $B$ then we say that
$f$ extends holomorphically along $L$ if $f\vert L\cap bB$ has a continuous extension through $L\cap\overline B$
which is holomorphic on $L\cap B$. Given $a\in\C^2$ we denote by $\L (a)$ the family of all complex
lines passing through $a$.
\vskip 2mm
\noindent\bf QUESTION 1.1\ \it
Let $a, b\in B, \ a\not=b$. Assume that $f\in C(bB)$ extends holomorphically along every complex
line in $\L (a) \cup \L (b)$. Must $f$
extend holomorphically through $B$? \rm 
\vskip 2mm
\noindent\bf EXAMPLE 1.1\ \ \rm Let $k\in\N$. For each $(z,w)\in bB$ let
$$
f(z,w) = \left\{\eqalign{ &z^{k+2}/\overline z\hbox{\ if\ }z\not= 0\cr
&0 \hbox{\ if\ } z=0 .\cr}\right.
$$
The function $f$ is of class $C^k$ on $bB$. Let $L$ be any complex line that meets the
disc $\{0\}\times\D$.
Then, if $L$ is not equal to the $w$-axis, we have $L\cap bB=\{(p_1+\z q_1, p_2+\z q_2)\colon\ \z\in b\D\} $
where $|p_1|<|q_1|$
so that the circle $\{p_1+\z q_1\colon\ \z\in b\D\}$ surrounds the origin. It is easy to see
that in this case the function $\z\mapsto
(p_1+\z q_1)^{k+2}/(\overline{p_1+\z q_1})\ \ (\z\in b\D)$ extends holomorphically through $\D$. Thus, $f$ extends
holomorphically along each complex line that meets $\{0\}\times\D$, yet $f$ does not extend
holomorphically through $B$. Thus, for each
$k\in\N$ there is a function $f\in C^k(bB)$ which provides a negative answer to the question
above. In the present paper we show that there
are no such functions of class $C^\infty$:
\vskip 2mm
\noindent \bf THEOREM 1.1\ \it Let $a,b\in\overline B,\ a\not= b.$ If a function
$f\in C^\infty (bB)$ extends holomorphically along each
complex line in $\L (a) \cup\L (b)$ that meets $B$ then $f$ extends holomorphically through $B$. \rm
\vskip 2mm
\noindent In other words, for each pair of points $a, b\in \overline B,\ \L (a)\cup \L (b)$ is a
test family for holomorphic extendibility
for $C^\infty (bB)$.

In the present paper we consider more general pairs of points $a, b$. Given $a, b\in C^2,\ a\not= b$,
denote by $\Lambda (a,b)$ the complex
line passing through $a$ and $b$. Our main result is the following
\vskip 2mm
\noindent\bf THEOREM 1.2\ \it Let $a, b\in \C^2, a\not= b$.

(A) Suppose that one of the points $a, b$ is contained in $B$.

(A1) If $<a|b>\not= 1$ then $\L (a)\cup \L (b)$ is a test family for holomorphic extendibility for $C^\infty (bB)$.

(A2) If $<a|b> = 1$ then
$$
\left. \eqalign
{&\hbox{there is a function\ } f\in C^\infty (bB)\ \hbox{which extends holomorphically}\cr
&\hbox{along every complex line in\ }\L (a)\cup \L (b)\ \hbox{that meets\ } B,\cr
&\hbox{yet\ } f \hbox{\ does not
extend holomorphically through\ } B .\cr} \right\} \eqno (1.1)
$$

(B) Suppose now that both points $a, b$ are contained in $\C^2\setminus B$

(B1) If $\Lambda (a, b)$ meets $B$ then $\L (a)\cup \L (b)$ is a test family for
holomorphic extendibility for $C^\infty (bB)$

(B2) If $\Lambda (a,b) $ misses $\overline B$ then (1.1) holds. \rm
\vskip 2mm
\noindent\bf REMARK  \rm Note that Theorem 1.1 follows from Theorem 1.2. The
methods used in the present paper do not apply to the case when $\Lambda (a,b)$ is tangent to $bB$ and
it will remain an open question for which pairs of points $a, b$ of this sort $\L (a)\cup \L (b)$ is
a test family for holomorphic extendibility for $C^\infty (bB)$.
\vskip 2mm
If $\varphi $ is a continuous function on a circle $\Gamma$ then we say that $\varphi $ extends holomorphically
from $\Gamma $ if it extends holomorphically through the disc bounded by $\Gamma$.

Given $\alpha \in\D$ denote by $C_\alpha$ the family of the circles obtained as the images under the Moebius map
$z\mapsto (\alpha - z)/(1-\overline\alpha z)$ of all circles in $\DD$ centered at the origin. To prove Theorem 1.2
we will have to prove the following new, one variable result:
\vskip 2mm
\noindent\bf THEOREM 1.3\ \it Let $\alpha, \beta \in\D,\ \alpha\not=\beta$, and
let $\varphi $ be a continuous function on
$\DD$ which extends holomorphically from each circle in $C_\alpha\cup C_\beta $.
Then $\varphi $ is holomorphic on $\D$. \rm
\vskip 4mm
\bf 2.\ Simplifying the geometry by using automorphisms of $B$ \rm
\vskip 2mm
Let $a\in B$. Write $P_0=0$ and
$$
P_a(z)= {{<z|a>}\over{<a|a>}}a\ \ ((z\in\C^2)\ \hbox{if \ } a\not= 0,
$$
$$
s_a= \sqrt{1-|a|^2},\ Q_a=I-P_a.
$$
The map $\varphi_a$ defined by
$$
\varphi_a(z)={{a-P_a(z)-s_aQ_a(z)}\over{1-<z|a>}}\ \ (z\in\C^2,\ <z|a>\not=1).
\eqno (2.1)
$$
is an automorphism of $B$, a homeomorphism of $\overline B$ and a $C^\omega$-diffeomorphism
of $bB$. it is a fractional
linear map of $\C^2$ that maps $\{ z\in\C^2 \colon \ <z|a>\not= 1\}$ onto itself and satisfies
$\varphi_a(a)=0,\ \varphi_a(0)=a$,\ $\varphi\circ\varphi = id$. It maps complex lines to complex lines [Ru].
In particular, if $b\in\C^2$ and $<b|a>\not=1$ then it maps $\L (b)$ to $\L (\varphi_a(b)$ where
$\varphi_a(b) $ is a point in $\C^2$ and if $<b|a>=1$ then it maps $\L (b)$ to a family of
parallel complex lines. In the special case when
$a= (\lambda, 0)$ where $\lambda \in\D$, we have
$$
\varphi_a(z,w)=\Biggl(
   {{\lambda - z}\over{1-\overline\lambda z}}, -{{\sqrt{1-|\lambda|^2}}\over{1-\overline\lambda z}}w\Biggr).
   \eqno (2.2)
   $$
It is known that every automorphism of $B$ is of the form $U\circ\varphi _a$ for some $a\in B$ and for
some unitary map $U$ [Ru].

Let $\varphi$ be an automorphism of
$B$. If $f\in C^\infty (bB)$ then $f\circ\varphi ^{-1}\in C^\infty (bB)$ and $f$ extends holomorphically along
each complex line in $\L (a) $
that meets $B$ if and only if $f\circ\varphi^{-1}$ extends holomorphically
along each complex line $\varphi (L),\ L\in \L (a)$, that meets $B$,
that is, along each complex line in $\L (\varphi (a))$ that meets $B$. Since $f\circ\varphi^{-1}$extends
holomorphically
through $B$ if and only if $f$ does it follows that $\L (a)\cup \L (b)$ is a test family for holomorphy
for $C^\infty (bB)$ if and
only if $\L (\varphi (a))\cup \L (\varphi (b))$ is a test family for holomorphy for $C^\infty (bB)$. This
observation simplifies the geometry
by applying automorphisms before beginning the proof of Theorem 1.2.

Consider (A1). Let $a\in B$ and $<a|b>\not=1$. Then $\varphi _a(a) =0$
and $\varphi _a(b) $ is a point in $\C^2$. Hence, after a composition with a
unitary map, we may, with no loss of generality assume that
$a= (0,0)$ and $b=(t,0)$ where $t>0$. Now, consider (B1). Suppose that
$|a|\geq 1,\ |b|\geq 1$ and assume that $\Lambda (a,b)$ meets $B$.
Then we may, after composition by a unitary map, assume that
$\Lambda (a,b) = \{ (z,w)\in\ C^2\colon\ z=\lambda\}$ where $0\leq \lambda<1$.
The transform (2.2) maps $\Lambda (a,b)$ to the $w$-axis. So, with no loss of generality
assume that $a, b$ are both on one coordinate axis, say on $z$-axis, that is
$a =(\alpha, 0),\ b=(\beta, 0),\ \ \alpha\not=\beta$, where $|\alpha|\geq 1, \ |\beta|\geq 1$.
One possibility is $|\alpha|=|\beta|=1$.
If one of $\alpha, \beta $, say $\alpha$, satisfies $|\alpha |>1$ then the map (2.2)
with $\lambda = 1/\overline\alpha $ maps
the complex lines through $(\alpha ,0)$ to the complex lines parallel to the $z$-axis
and the point $(\beta ,0)$ to a point on
$z$-axis, which, after a a suitable rotation $(z,w)\mapsto (e^{i\omega}z, w)$ we may
assume that is of the form $(t,0)$ where $t\geq 1$. Thus, (A1) and (B1) will be proved once we
have proved that each of the following
families of lines (2.3), (2.4) and (2.5) is a test family for holomorphic extendibility for $C^\infty (bB)$:
$$
\left. \eqalign{
   &\hbox{the complex lines passing through the origin and the complex}\cr &\hbox{lines passing through
   a point}\ (t,0) \hbox{\ where\ } t>0, \cr}\right\}
   \eqno (2.3)
$$
$$
\left. \eqalign{
   &\hbox{the complex lines passing through the point\ } (\alpha, 0) \hbox{\ and the
   complex}\cr &\hbox{lines passing through the point}\ (\beta ,0) \hbox{\ where\ }
   |\alpha|=|\beta|=1,\ \alpha\not=\beta,\cr}\right\}
   \eqno (2.4)
$$
$$
\left. \eqalign{
   &\hbox{the complex lines parallel to the } z-\hbox{axis\ } \hbox{and the complex}\cr &\hbox{lines
   passing through a point}\ (t,0) \hbox{\ where\ } t\geq 1 .\cr}\right\}
   \eqno (2.5)
$$

Consider now (A2), so let $a\in B$ and $<a|b>=1$. The map $\varphi_a$ maps $\L (a)$ to $\L (0)$
and $\L (b)$ to a family of parallel lines for which, after a composition with a unitary map, we may assume that
are parallel to the $z$-axis. Now, consider (B2). After composition by a unitary map we may assume that $\Lambda (a,b) =
\{(z,w)\in\C^2\colon z=\mu\}$ where $\mu >1$  so that $a=(\mu,\alpha)$, \ $(b=(\mu,\beta)$. An easy computation shows that 
the map (2.2) with
$\lambda = 1/\mu$ maps $\L (a)$ to the complex lines parallel to $(1,\eta_\alpha)$ and $\L (b)$  to the
complex lines parallel to $(1,\eta_\beta)$ where
$$
\eta_\alpha = {{\alpha}\over{\sqrt{\mu^2-1}}},\ \ \eta_\beta = {{\beta}\over{\sqrt{\mu^2-1}}}.
$$
It follows that (A2) will be proved once we have proved that
$$
\left. \eqalign{
&\hbox{there is an \ } f\in C^\infty (bB)\hbox{\ that
extends holomorphically along every complex}\cr
&\hbox{line passing through the origin and every
complex line parallel to}\ z-\hbox{axis}\cr
&\hbox{which meets }\ B\hbox{\ yet\ } f\hbox{\ does
not extend holomorphically through}\ B,\cr}\right\}\eqno (2.6)
$$
and (B2) will be proved once we have proved that
$$
\left. \eqalign
{&\hbox{given\ }p, q\in\C,\ p\not= q, \ \hbox{there is an\ } f\in C^\infty (bB)\hbox{\ which extends
holomorphically }\cr
&\hbox{along every complex line parallel to\ } (1,p)\ \hbox{which meets\ }B,\
\hbox{and along every}\cr
&\hbox{complex line parallel to
\ }(1,q) \hbox{\ which meets\ }
B,\ \hbox{yet \ }f\hbox{\ does not\ }
\hbox{extend}\cr &\hbox{holomorphically through\ }B.\cr}\right\}\eqno (2.7)
$$
It is easy to see that the function $f(z,w)=|w|^2$ satisfies (2.6). To get an example
of a function $f$ that satisfies (2.7) we follow A.\ M.\ Kytmanov and S.\ G.\ Myslivets [KM] and put
$$
f(z,w) = \overline z[z(1+|p|^2)+\overline p(w-pz)].[z(1+|q|^2)+\overline q(w-qz)].
$$
It is easy to check that $f$ satisfies (2.7): If a complex line is parallel to $(1,p)$ then
it has the form $\{(\z, c+\z p)\colon \ \z\in\C\}$ for some $c\in\C$. Note that $|\z|^2+|c+\z p|^2 = 1 $ implies that
$\overline\z [\z(|p|^2+1)+\z\overline p]=1-|c|^2-\z p\overline c $ so for such $\z$, \ $f(\z ,c+\z p)=
[1-|c|^2-\z p\overline c].[\z (1+|q|^2)+\overline q (c+\z p-\z q)]$ depends holomorphically on $\z $.
Repeating the reasoning for $(1,q)$ we see that $f$ satisfies (2.7).

Note that both examples above are real analytic.

It remains to prove that the families (2.3), (2.4) and (2.5) are test families for holomorphic extendibility for $C^\infty (bB)$.
\vskip 4mm
\bf 3.\ Reduction to a sequence of one variable problems \rm
\vskip 2mm
As in [G1, G2] we shall use the Fourier series decomposition and averaging. Suppose that $f\in C^\infty (bB)$. Given $n\in Z$
and $z\in \D$ let
$$
c_n(z)= \left( {1\over{\sqrt{1-|z|^2}}}\right)^n{1\over{2\pi}}
\int_{-\pi}^{\pi}e^{-in\theta}f(z,e^{i\theta}\sqrt{1-|z|^2}) d\theta .
$$
so that, since $e^{i\theta}\mapsto f(z,\sqrt{1-|z|^2}e^{i\theta})$ is smooth, we have
$$
f(z,e^{i\theta}\sqrt{1-|z|^2}) = \sum_{n=-\infty}^\infty (\sqrt{1-|z|^2})^n c_n(z)e^{in\theta} = \sum_{n=-\infty}^\infty
c_n(z)(\sqrt{1-|z|^2}e^{i\theta})^n .
$$
If $(z,w)\in bB$ then writing $w=e^{i\theta}\sqrt{1-|z|^2}$ the preceding discussion implies that
$$
f(z,w)=\sum_{n=-\infty}^\infty c_n(z)w^n\ \ ((z,w)\in bB,\ w\not=0).
$$
The coefficients  $c_n$ are continuous on $\D$ and from the definition it follows that if $n\leq 0$ they
also continuously extend to $\overline D$. We shall show that when $f\in C^\infty (bB)$ the same holds for $n>0$:
\vskip 2mm
\noindent\bf LEMMA 3.1\ \it Suppose that $f\in C^\infty (bB)$. Then for each
$n\in Z$ the function $z\mapsto c_n(z)\ (z\in\D )$ extends continuously to $\DD$. \rm
\vskip 2mm
Let $z_0\in\C$ and assume that $f\in C(bB)$ extends holomorphically along every complex
line passing through $(z_0,0)$ which meets $B$. If $L$ is such a
a complex line then so is $\omega_\theta (L)$ where, for $\theta\in\R$,
\ $\omega_\theta (z,w)=(z,e^{-i\theta} w)$. It follows that $f$ extends holomorphically along every
complex line $\omega_\theta (L)$ which is the same to say that for each
$\theta\in\R,\ (z,w)\mapsto f(z,e^{i\theta}w)$ extends holomorphically along $L$ and hence for each $n\in Z$,
the same holds for 
$$
\Psi_n(z,w)={1\over{2\pi}}\int_{-\pi}^\pi e^{-i n\theta}f(z,we^{i\theta})d\theta,
$$
a continuous function on $bB$. Note that
$$
\Psi _n(z,w)= w^nc_n(z)\ \ ((z,w)\in bB,\ w\not=0).
$$
This proves
\vskip 2mm
\noindent\bf PROPOSITION 3.1\it \ If a function $f\in C(bB)$ extends holomorphically along each
complex line that passes through $(z_0,0)$and meets $B$ then for each $n\in Z$ the same
holds for the continuous extension of $(z,w)
\mapsto w^n c_n(z)\ ((z,w)\in bB,\ w\not= 0)$ to $bB$.
\vskip 2mm
\noindent \rm Analogous statement holds for complex lines parallel to $z-$axis. For such
complex lines $w=const$ so we have
\vskip 2mm
\noindent\bf PROPOSITION 3.2\ \it If a function $f\in C(bB)$ extends holomorphically along
each complex line that is parallel to
the $z-$axis and meets $B$ then for each $n\in Z$ the same holds for the function
$(z,w)\mapsto c_n(z)\ ((z,w)\in bB,\ w\not=0)$.
\vskip 2mm
\noindent \rm Note that this is the same as to say that $c_n$ extends holomorphically from
each circle in $\D$ centered at the origin.

Denote by $\pi_1$ the projection onto $z-$axis, \ $\pi_1(z,w)=z$.

Assume now that for given $z_0\in\C$ a function $f\in C(bD)$ extends holomorphically
along each complex line passing through $(z_0,0)$.
By Proposition 3.1 we know that for each $n\in Z$ the continous extension to $bB$ of the function
$(z,w)\mapsto w^nc_n(z)$ has the same property. So let $L$ be such a complex line that is not the $z$-axis 
and is not parallel to the $w$-axis so that $\pi_1(L\cap bB)=\{ p+\z q\colon\ \z\in\bD\}$ 
is not a point, and let $n\in Z$. Write
$L\cap bB =\{ (p+\z q,r+\z s)\colon\ \z\in b\D\} $, note that $q\not=0,\ s\not=0$ and note that $c_n(p+\z q)$ is not
defined for at most one point $\z\in \bD$ for which
we have $r+\z s=0$ and that this can happen only in the case when $|z_0|=1$ when $p+\z q=z_0$.  Thus,
$$\left.\eqalign{&\hbox{the continuous extension of\ }
\z\mapsto (r+\z s)^nc_n(p+\z q) \hbox{\ to\ }\bD\cr
&\hbox{extends holomorphically through }\ \D.\cr}\right\} \eqno (3.1)
$$
Since $L$ passes through $(z_0,0)$ there is a $\z _0\in\C$ such that $p+\z_0 q=z_0$ and $r+\z _0 s=0$.
Clearly $\z_0= (z_0-p)/q$ so that $r+(z_0-p)s/q = 0$. Writing $p+\z q=z$
we get $\z = (z-p)/q$ and it follows that  $r+\z s=(z-z_0)s/q$,  so (3.1) implies that the
continous extension of
$z\mapsto (z-z_0)^nc_n(z)$ to the circle $\{p+\z q\colon\ \z\in b\D\}$\ extends
holomorphically from this circle. This proves
\vskip 2mm
\noindent\bf PROPOSITION 3.3\ \it If $f\in C(bB)$ extends holomorphically along each
complex line passing through $(z_0,0)$ that meets $B$ then for each such
complex line $L$ and for each $n\in Z$ the function $z\mapsto (z-z_0)^nc_n(z)$ extends
holomorphically from $\pi_1(L\cap bB)$.
\vskip 2mm \rm
\noindent\bf REMARK \rm In the case when $z_0\in \bD$  the circle $\pi_1(L\cap bB)$
contains $z_0$ and in this case we have to be more precise and say that
the continous extension of $z\mapsto (z-z_0)^n c_n (z)$
to the circle $\pi_1(L\cap bB)$ extends holomorphically from this circle.
\vskip 2mm
We shall prove the following
\vskip 2mm
\noindent\bf LEMMA 3.2\ \it Suppose that $\varphi $ is a continuous function on $\DD$.
Suppose that $\F $ is one of the families (2.3), (2.4), (2.5) of complex lines. Suppose that
for each complex line belonging to $\F $ which meets the ball,\  $\varphi $ extends holomorphically from the
circle $\pi_1(L\cap bB)$. Then $\varphi $ is holomorphic on $\D$. \rm
\vskip 2mm
Assuming Lemma 3.1 and Lemma 3.2 for a moment we can now complete the proof
of Theorem 1.2 as follows:

Let $f\in C^\infty (bB)$. By Lemma 3.1 for each $n\in Z$ the function $c_n$
extends continuously to $\DD$.

Suppose that $f$ extends holomorphically along each complex line belonging to one of
the families (2.3), (2.4), (2.5) that meets $B$. Suppose that $n\geq 0$. By Proposition 3.3 and Proposition 3.2 
either
$$
\left.\eqalign{
   &\hbox{for each \ }L\hbox{\ belonging to the family (2.3) the function}\
   z\mapsto \varphi (z) = \cr
   &z^n(z-t)^nc_n(z)
   \hbox{\ extends holomorphically from\ } \pi_1(L\cap bB)\cr}\right \} \eqno(3.2)
   $$
or
$$
\left.\eqalign{
   &\hbox{for each \ }L\hbox{\ belonging to the family (2.4) the function}\
   z\mapsto \varphi (z) = \cr
   &(z-\alpha)^n(z-\beta)^nc_n(z)
   \hbox{\ extends holomorphically from\ } \pi_1(L\cap bB)\cr}\right \} \eqno(3.3)
   $$
or
$$
\left.\eqalign{
   &\hbox{for each \ }L\hbox{\ belonging to the family (2.5) the function}\
   z\mapsto \varphi (z) = \cr
   &(z-t)^nc_n(z)
   \hbox{\ extends holomorphically from\ } \pi_1(L\cap bB . )\cr}\right \} \eqno(3.4)
   $$
   
   \noindent Here we used the fact that if a function $\psi$ extends holomorphically from a 
   circle $\Gamma$ and if $p\in\C$ then the function $z\mapsto (z-p)^n\psi (z)$ extends holomorphically from $\Gamma$. 
   Since in each of the three cases the function $\varphi$ is continous on $\overline\D$ it follows by Lemma 3.2 that
    in each of the three cases the function $\varphi$ is holomorphic on $\D$. It follows that $c_n$ is holomorphic on $\D$ 
    except perhaps at two poles, $0$ and $t$, which, by the continuity of $c_n$, are removable singularities, hence 
    $c_n$ is holomorphic on $\D$.
    
    Now let $n<0$. Then, if for a function $\psi$ the function $z\mapsto (z-p)^n\psi (z)$ extends holomorphically 
    from a circle $\Gamma $ it follows that the function $\psi $ extends holomorphically 
    from $\Gamma $.  Together with Proposition 3.3 and Proposition 3.2 
    this implies that the function $c_n$ extends holomorphically 
    from $\pi_1(L\cap bB)$ either for each $L$ belonging to the family (2.3), 
    or for each $L$ belonging to the family (2.4) or for each $L$ belonging to the family (2.5).
    Since $c_n$ is continuous on $\D$ it follows by Lemma 3.2 that
    in each of the three cases the function $c_n$ is holomorphic on $\D$.
    
It follows that
$$
f(z,w) = \sum_{n=0}^\infty w^nc_n(z) \ \  ((z,w)\in bB)
$$
where each $c_n$ is continuous on $\DD$ and holomorphic on $\D$. Since this
is a Fourier decomposition of a
smooth function, the series converges uniformly on $bB$ and so by the maximum
principle, uniformly on $\overline B$. So $f$ extends holomorphically
through $B$. This completes the proof of Theorem 1.2. It remains to prove
Lemma 3.1 and Lemma 3.2.
\vskip 4mm
\bf 4.\ Proof of Lemma 3.1\rm
\vskip 2mm
Assume that $f\in C^\infty (bB)$. We have to show that for each $n\in\N$ the function
$$
z\mapsto c_n(z)= \left\{ {1\over{\sqrt{1-|z|^2}}}\right)^n{1\over{2\pi}}
\int_{-\pi}^{\pi}e^{-in\theta}f(z,e^{i\theta}\sqrt{1-|z|^2}) d\theta ,
$$ which is defined and smooth on $\D$, extends continuously to $\DD$.

Write $z=\sqrt{1-R^2}e^{i\varphi}$ We have $\sqrt{1-|z|^2} = R$ so
$$
c_n(z) =
{1\over{R^n}}{1\over{2\pi}}
\int_{-\pi}^{\pi}e^{-in\theta}f(\sqrt{1-R^2}e^{i\varphi},Re^{i\theta}) d\theta .
\eqno (4.1)
$$
We will prove the continuous extendibility of $c_n$ to $\DD$ by showing that as $R\searrow 0$,
the functions $e^{i\varphi}\mapsto c_n(\sqrt{1-R^2}e^{i\varphi})$ converge to a function
$e^{i\varphi}\mapsto c_n(e^{i\varphi})$, uniformly in $\varphi,\ -\pi\leq \varphi\leq \pi$.
With no loss of generality assume that $f\in C^\infty (C^2)$. To prove the continuous
extendibility notice first that whenever $\Phi $
is a smooth function on $\C^2$ and $j\in\N$ then the integration by parts gives
$$
\eqalign{
   & \int_{-\pi}^{\pi}\Phi (\sqrt{1-R^2}e^{i\varphi}, Re^{i\theta})e^{-ij\theta}d\theta =\cr
   &= {1\over{ij}}\int_{-\pi}^{\pi}{\partial\over{\partial\theta}}
   \left[\Phi(\sqrt{1-R^2}e^{i\varphi},Re^{i\theta}
   \right]e^{-ij\theta}d\theta =\cr
   &= {1\over{ij}}\int_{-\pi}^\pi \Bigl[
   {{\partial\Phi}\over{\partial w}}(\sqrt{1-R^2}e^{i\varphi},
    Re^{i\theta})iRe^{i\theta} -\cr
   &-
   {{\partial\Phi}\over{\partial {\overline w}}}(\sqrt{1-R^2}e^{i\varphi}, Re^{i\theta})
   (-iRe^{-i\theta}) \Bigr] e^{-ij\theta}d\theta =\cr
   &= {R\over j}\int_{-\pi}^\pi{{\partial\Phi}\over{\partial w}}(\sqrt{1-R^2}e^{i\varphi},
   Re^{i\theta
   })e^{-i(j-1)\theta}d\theta -\cr
   &- {R\over j}\int_{-\pi}^\pi{{\partial\Phi}\over{\partial {\overline w}}}(\sqrt{1-R^2}
   e^{i\varphi},
   Re^{i\theta})e^{-i(j+1)\theta}d\theta \cr}
   $$
   where both integrands have the same form as the integrand at the beginning.

   Given $n\in\N$ we apply the preceding reasoning $n$ times to see that
   $$
   R^nc_n(z)={1\over{2\pi}}\int_{-\pi}^\pi e^{-in\theta}f(\sqrt{1-R^2}e^{i\varphi}, Re^{i\theta})d\theta
   $$
   is a finite sum of terms of the form
   $$R^n\gamma_{j,k}\int_{-\pi}^\pi{{\partial^nf}\over{\partial w^j\partial\overline w}^k}(\sqrt{1-R^2}e^{i\varphi},
   Re^{i\theta})e^{-in\theta}e^{ij\theta}e^{-ik\theta}d\theta
   $$
   where
   $j+k=n$ and where where $\gamma_{j,k}$ are constants, so $c_n(z)$ is a finite sum of terms of the form
   $$
   \gamma_{j,k}\int_{-\pi}^\pi{{\partial^nf}\over{\partial w^j\partial\overline w}^k}(\sqrt{1-R^2}e^{i\varphi},
   Re^{i\theta})e^{-i(n-j+k)\theta} d\theta .
   $$
   As $R\searrow 0$ each of these integrals converges, uniformly in $
   \varphi, \ -\pi\leq \varphi\leq\pi$, to

   $$
   \gamma_{j,k}\int_{-\pi}^\pi{{\partial^nf}\over{\partial w^j\partial\overline w^k}}
   (\sqrt{1-R^2}e^{i\varphi},
   0)e^{-i(n-j+k)\theta} d\theta
   $$
   which equals $0$ if $ n-j+k>0 $ and
   $$
   \gamma_{n,0}\int_{-\pi}^\pi{{\partial ^nf}\over{\partial w^n}}(e^{i\theta},0)d\theta
   $$ if $n-j+k=0$ which happens when $j=n,\ k=0$.
 Thus, as $R\searrow 0,\ c_n(\sqrt{1-R^2}e^{i\theta})$ converges, uniformly for $-\pi\leq \varphi\leq \pi$ to
 $$
 2\pi \gamma_{n,0}{{\partial ^nf}\over{\partial w^n}}(e^{i\varphi}, 0).
 $$
 This completes the proof.
 \vskip 4mm
 \bf 5.\ Projections of the intersection of bB with complex lines \rm
 \vskip 2mm
 In this section we give a precise description of the circles $\pi_1(L\cap bB)$ where $L$ is a complex line passing
 through $(t,0)$ where $t\geq 0$.
 Given $z\in\C$ and $r>0$ write $\D (z,r) = \{\z\in\C\colon\ |\z-z|<r\}$.

 The intersections of $bB$ with complex lines through the origin are
 $$
 \{(R\z,e^{i\omega}\sqrt{1-R^2}\z)\colon\ \z\in b\D\},\ 0\leq R\leq 1,\ \omega\in\R
 \eqno (5.1)
 $$
 and the intersections of $bB$ with the complex lines through $(t,0),\  0<t<1$, are the images of (5.1) under the Moebius map
 $$
 \varphi (z,w)=\left({{t-z}\over{1-tz}}, -\sqrt{1-t^2}{w\over{1-tz}}\right)
 $$
 which takes the origin to the point $(t,0)$. Fix $t,\ 0<t<1$. Then $\pi_1 (L\cap bB)$ for complex lines passing through
 $(t,0)$ are the circles
 $$
 \left\{ {{t-R\z}\over{1-tR\z}},\ \z\in\bD\right\} = 
 \left\{ {{R\z+t}\over{1+tR\z}},\ \z\in\bD\right\},\ 0 < R\leq 1.
 $$
 Given $R$, the diameter of such a circle is the closed interval $[(-R+t)/(1-tR), (R+t)/(1+tR)]$ on the real axis so the center is
 $$
 T= {1\over 2}\left( {{-R+t}\over{1-tR}} + {{R+t}\over{1+tR}}\right) = {{t(1-R^2)}\over{1-t^2R^2}}
 $$
 and the radius is
 $$
 \rho = {1\over 2}\left({{R+t}\over{1+tR}} - {{-R+t}\over{1-tR}}\right) = {{R(1-t^2)}\over{1-t^2R^2}}.
 $$
 Notice that when $R$ increases from $0$ to $1$,\ $T$ decreases from $t$ to $0$. Since $R^2=(T-t)/(t(Tt-1))$ it
 follows that $\rho^2=(T-t)(T-1/t)$. This shows that if $0<t<1$ then for the complex lines $L$ passing through the point
 $(t,0)$, the circles $\pi_1(L\cap bB)$ are
 $$
 \bD (T,\sqrt\Tt),\ \ 0\leq T\leq t.
 \eqno(5.2)
 $$
If $t=1$ then for complex lines $L$ passing through the point $(1,0)$, \ $\pi_1(L\cap bB)$
are the circles
contained in $\DD$ which pass through the point $1$.

If $t>1$ then for complex lines $L$, passing through the point $(t,0)$ which meet $B$, the
circles $\pi_1(L\cap B))$ are, similarly to (5.2), the circles
$$
\bD (T,\sqrt\Tt),\ 0\leq T<1/t.
\eqno (5.3)
$$
To verify this, one writes $z=x+iy,\ w=u+iv$ and looks at the $(x,u)-$plane
$E=\{ (x,u)\in\C^2\colon\ x\in\R, y\in\R\}$. Let $\ell$ be a
line in $E$ passing through the point $(t,0),\ t>1$, which intersects the
open unit disc in $E$; denote the intersection by $J$. The segment
$J$ is the diameter of the disc $L\cap B$ where $L$ is the complex line in $\C^2$
that contains $\ell$ and the projection
of $J$ in $E$ to the $x-$axis, the real axis in the $z-$axis, is the diameter
of $\pi_1(L\cap B)$. If $T$ is its midpoint then a simple calculation
in $\R^2$ shows that the length of $\pi_1(J)$ is $2\sqrt\Tt$.

This shows that to prove Lemma 3.2 we have to consider three different families
(actually pairs of families) of circles in $\DD$:
\vskip 1mm
\noindent\it $(\cC_1)$\ for $\alpha, \beta \in \bD,\ \alpha\not=\beta$, the
family of all circles in $\DD$ passing through $\alpha$
and the family of all circles in $\DD$ passing through $\beta$
\vskip 1mm
\noindent $(\cC_2)$\ the family of all circles
centered at the origin and the family of all circles in $\DD$ passing through 1
\vskip 1mm
\noindent $(\cC_3)$\ the family of all circles centered at the origin and, for $t,\ 0<t<1$, the family
$$
b\D(T,\sqrt \Tt ),\ 0\leq T<t,
$$
that is, the family of all circles obtained from the circles centered at the
origin by the Moebius transform $z\mapsto (t-z)/(1-tz)$.
\vskip 1mm \rm

Lemma 3.2 will be proved once we have proved, for each $j,\ 1\leq j\leq 3$,
that if $\varphi $ is a continuous function on $\DD$ which extends
holomorphically from each circle belonging to $\cC_j$, then $\varphi$ is
holomorphic on $\D$. For $\cC_2$ this is the main result of [G6]. As
mentioned in [G6] it is proved in the same way for $\cC_1$. It remains to prove it for $\cC_3$:
\vskip 2mm
\noindent\bf LEMMA 5.1\ \it Let $0<t<1$. Suppose that $\varphi\in C(\DD )$ extends
holomorphically from each circle in $\DD$ centered at
the origin and from each circle $\bD (T,\sqrt\Tt),\ 0\leq T<t$. Then $\varphi$ is
holomorphic on $\D$. \rm
\vskip 2mm
\noindent It is clear that by proving Lemma 5.1 we also prove Theorem 1.3. In fact,
the statements of Lemma 5.1 and Theorem 1.3
are equivalent.

Lemma 3.2 will be proved once we have proved Lemma 5.1. To prove Lemma 5.1 we will
first use semiquadrics as in [AG,G3]
to formulate the problem in $\C^2$ and then show that the idea of A.Tumanov [T2] to
use an argument of H.Lewy [L] together with the Liouville theorem
still aplies in our situation. The proof will be similar to the proof of the main
result of [G6] but more complicated.
\vskip 4mm
\bf 6.\ Proof of Lemma 5.1, Part 1 \rm
\vskip 2mm
We introduce semiquadrics
to pass to an associated problem in $\C^2$. Given $a\in\C$ and $r>0$ let
$$
\Lambda (a,r) = \{ (z,w)\in\C^2\colon\ (z-a)(w-\overline a)=r^2,\ 0<|z-a|<r\} .
$$
be the semiquadric associated with the circle $b\D (a,r)$. Write
$\Sigma = \{ (\z, \overline\z )\colon\ \z\in\C\}$. \ $\Lambda(a,r)$  is a closed complex submanifold
of $\C^2\setminus\Sigma $ which is attached to $\Sigma $ along
$b\Lambda (a,r)= \{ (\z,\overline \z)\colon\ \z\in\bD (a,r)\}$. A
continuous function $g$ extends holomorphically from the circle $\bD (a,r)$ if and
only if the
function $G$, defined on $b\Lambda (a,r)$ by
$G(\z,\overline\z )= g(\z )\ (\z\in b\D (a,r))$ has a bounded continuous extension
to $\overline {\Lambda (a,r)} =
\Lambda (a,r)\cup b\Lambda (a,r)$ which is holomorphic on $\Lambda (a,r)$. In fact,
if we denote by the same letter $g$
the holomorphic extension of $g$ through $\D (a,r)$ we have
$$
G\Bigl(z,\overline a + {{r^2}\over{z-a}}\Bigr) = g(z)\ \ (z\in\DD(a,r)\setminus \{ a\})
$$
and  if we define $G(a,\infty)=g(a)$ we get a continuous function $G$ on
$\overline {\Lambda (a,r)}\cup \{(a,\infty)\}$, the closure of $\Lambda (a,r)$ in
$\C\times \overline\C$. It is known that if $(a,r)\not= (b,\rho)$ then $\Lambda (a,r)$
meets $\Lambda (b,\rho)$ if and only if $a\not= b$ and one of the circles
$b\D (a,r),\ b\D(b,\rho)$ surrounds the other. If this happens then $\Lambda(a,r)$
and $\Lambda (b,\rho)$ meet at precisely one point [G3].

We begin with the proof of Lemma 5.1. Let $\varphi $ and $t$ be as in Lemma 5.1.
By our assumption, $\varphi $ extends holomorphically from two
families of circles:\ $\{ \bD(0,R)\colon\ 0<R\leq 1\}$ and $\{ \bD (T, \sqrt\Tt ),\ 0\leq T<t\}$. Accordingly,
there are two families of semiquadrics:
\ $\{ \Lambda (0,R), 0<R\leq 1\}$ and $\{ \Lambda (T,\sqrt\Tt )\colon\ 0\leq T<t\}$,\
and the function $\Phi (\z ,\overline\z ) =
\varphi (\z )$, defined on
$\{(\z,\overline\z )\colon\ \z\in\DD\}$ has a bounded holomorphic extension through
each of these semiquadrics.

Consider the first family. In this family the semiquadrics are pairwise disjoint.
Let $L$ be the closure of their union in $[\C\setminus\{ 0\}]\times \C$:
$$
L = \bigcup_{0<r\leq 1} [\Lambda (0,R)\cup b\Lambda (0,R)].
$$
The continuity of $\varphi $ together with the maximum principle implies that our
function $\Phi $ extends from
$L\cap\Sigma=\{ (\z,\overline\z)\colon\ \z\in\DD\setminus\{ 0\}\} $ to a bounded
continuous function $\Phi $ on $L$ so that
the extension $\Phi $ is holomorphic on each fiber $\Lambda (0,R)$. Note that $L$
is a CR-manifold in $[\C\setminus\{ 0\}]\times\C$
with piecewise smooth boundary consisting of two pieces, $\Lambda (0,1)$ and
$\{ (\z,\overline\z )\colon\ \z\in\DD\setminus\{ 0 \}\} $ and the function $\Phi$,
being holomorphic on fibers, is CR on its interior $L_0=\cup_{0<R<1}\Lambda (0,R)$, that is
$$\int_{L_0}\Phi\overline\partial\omega = 0
$$
for each smooth, (2,0)-form $\omega $ on $\C^2$ whose support meets $L_0$ in a compact set.

Now, look at the second family of semiquadrics
$$
\eqalign{
\Lambda \left(T,\sqrt\Tt \right) = \{(z,w)\in\C^2\colon\ (z-T)(w-T)=\Tt, &\cr
0<|z-T|<\sqrt\Tt\}, \ 0\leq T<t. &\cr}
$$
Observe that these semiquadrics are not pairwise disjoint since they all contain
the point $(t,1/t)$. Since two semiquadrics can meet at
at most one point it follows that the sets
$$
\Lambda \left(T,\sqrt\Tt\right)\setminus\{ (t,1/t)\} ,\  0<T<t,
$$
are pairwise disjoint and so their union
$$
N_0 =\left[\bigcup_{0<T<t}\Lambda (T,\sqrt\Tt)\right] \setminus\{ (t,1/t)\}
$$ is a CR manifold. Let $N$ be the closure of $N_0$ in $[\C\setminus\{t\}]\times\overline\C$, that is,
$$
N = \bigcup_{0\leq T < t}\left[\Lambda (T,\sqrt\Tt)\cup b\Lambda (T,\sqrt\Tt)\cup\{ (T,\infty)\}\right].
$$
Again, the continuity of $\varphi$ together with the maximum principle implies that
our function $\Phi$
extends from $N\cap\Sigma = \{(\z,\overline\z)\colon\ \z\in\DD\setminus\{ t\} \}$
continuously to $N$ so that our extension $\Phi$ is holomorphic on each fiber $\Lambda (T,\sqrt\Tt),\
0\leq T<t$. The part of $N$ contained in $\C\times\C$ is a smooth CR manifold in
$[\C\setminus\{ t\} ]\times \C$
 with boundary consisting of two pieces, $\Lambda (0,1)$ and $N\cap\Sigma $ and
 the function $\Phi $ is CR in the interior $N_0$.
 \vskip 4mm
 \bf 7.\ Proof of Lemma 5.1, Part 2\rm
 \vskip 2mm
 Let $\varphi$ and $t$ be as in Lemma 5.1. We have shown that the function $\Phi $ extends continuously from
 $\{ (\z,\overline\z)\colon\ \z\in\D\setminus \{ 0\}\}$ to $L$ and from
 $\{ (\z,\overline\z)\colon\ \z\in\D\setminus \{ t\}\}$
 to $N$ so that the extensions are holomorphic on fibers of $L$ and $N$.
 We would like that these extensions define a function $\Phi$ on $L\cup N$.
 However, this is not possible since $L$ and $N$ intersect. There is one
 piece of $L\cap N$,
 namely $\Lambda (0,1)$ on which both extensions coincide. There are other
 points of $L\cap N$. Obviously, all such points
 of intersection are of the form $(x,y)$ where $x\in\R,\ y\in\R$.
 We now show that there is an $\eta >0$ such that there are no such intersections
 with $0<x<\eta$.

 For each $R,\ 0<R\leq 1$, write $\cT _R =\Lambda(0,R)$ and for each $T,\ 0\leq T<t$, write $\cS_T =
 \Lambda (T,\sqrt\Tt)$. Clearly $\cS_0=\cT_1$. Choose $\eta$ so that
 $$
 0<\eta<{1\over{2(t+1/t)}}
 \eqno (7.1)
 $$
 and notice that $\eta<1/2$ and $\eta<t$.
 \vskip 2mm
 \noindent\bf PROPOSITION 7.1\ \it Let $\eta$ satisfy \rm(7.1)\it and assume that
 $0<x<\eta$. Then $\{ x \}\times\C$ contains no point of
 $\cS_T\cap\cT_R$ with $0<T<t$ and $0<R<1$.
 \vskip 1mm
 \bf Proof.\ \rm We shall prove the proposition by proving that
 $$
 \hbox{if\ } 0<R<1\hbox{\ and if\ }(x,y)\in\cT_R\hbox{\ then\ }x<y<1/x,
 \eqno (7.2)
 $$
 $$
 \hbox{if\ }0<T<t\ \hbox{and if\ } (x,y)\in\cS_T\hbox{\ then either\ } 1/x<y<\infty\hbox{\ or\ }-\infty <y<x.
 \eqno (7.3)
 $$
 A simple picture shows that there is no $(x,y)\in\cT_R$ if $R\leq x$
 and when $x<R<1$ than there is precisely one $y$ such
 that $(x,y)\in\cT_R$ and
 $$
 \hbox{as\ }R\hbox{\ moves from\ }x{\ to \ }1,\ y\hbox{\ moves from\ }x\hbox{\ to\ }1/x.
 \eqno (7.4).
 $$
 This takes care of (7.2). We now turn to (7.3). We first determine the
 interval $(0,T_0)$ of those
 $T$ for which $x$ is contained in the disc $\D (T,\sqrt\Tt )$, that is,
 of all those $T$ for which there is a point
 $(x,y)
 \in\cS_T\cup\{ (T,\infty )\} $. Clearly $T_0$ is determined by the condition
 that $(x,x)\in b\cS_{T_0}$, that is,
  $(x-T_0)^2= (T_0-t)(T_0-1/t)$ which gives
  $$
  T_0={{1-x^2}\over{(t+1/t)-2x}}.
  $$
  Since $x<\eta <1/2$, \ (7.1) implies that $x(t+1/t)<1/2$ hence
  $$
  x^2-x(t+1/t)+1>0.
  \eqno (7.5)
  $$
  Obviously $x<T_0$. One can see this also by rearranging (7.5) to $1-x^2>x(t+1/t)-2x^2$
  which gives $T_0>x$.

  Now, for each $T,\ 0\leq T\leq T_0$ we compute the unique $y=y(T)$
  such that $(x,y)\in\cS_T\cup\{ (T,\infty )\}$.
  It
  is clear that $y(0) = 1/x$ and $y(x) = \pm\infty$. We now show that
  $$
  \left.\eqalign{
     &\hbox{as\ }T\hbox{\ increases from\ }0\hbox{\ to\ }x,\ y(T) \hbox{\ increases from\ }1/x\hbox{\ to\ }+\infty \cr
     &\hbox{and as\ }T\ \hbox{increases from\ } x\hbox{\ to \ }T_0,\ y(T)\hbox{\ increases from\ }-\infty \hbox{\ to\ }x.
     \cr}\right\}
     \eqno (7.6)
     $$
 To show this it suffices to show that ${{dy}\over{dx}}>0 $ for all $T,\ 0<T<T_0,\ T\not= x$. We have
 $$
 y(T) = T+{{1-(t+1/t)T+T^2}\over{x-T}}
 $$
 and a short computation shows that
 $$
 {{dy}\over{dT}}={{x^2-(t+1/t)x+1}\over{(x-T)^2}}
 $$
 which, by (7.5) is always positive. This completes the proof of (7.3)
 and so completes the proof of Proposition 7.1.

 We have extended $\Phi$ to both $L$ and $N$ and so this common extension is
 well defined on $(L\cup N)\setminus
 (L\cap N)$. Thus, by Proposition 7.1, the extension is well defined on
 $$
 M=(L\cup N)\setminus \left[ ([-1,0]\cup [\eta,1])\times\overline C\right].
 $$
 It is continuous and holomorphic on fibers, that is, it is CR in the interior
 of the part of $M$ in $\C^2$. It is on this $M$ where
 we will apply the idea of Tumanov.
 \vskip 4mm
 \bf 8.\ Completion of the proof of Lemma 5.1\rm
 \vskip 2mm
 Let $\cS =\D\setminus \left( (-1, 0]\cup [\eta, 1)\right)$. Given $z\in\cS$, let $M_z=\{ \z\in\C\colon\ (z,\z )\in M\}$.

 Given $z\in\cS,\ z\not\in\R$, let $C_z$ be the circle passing through
 $t, 1/t$ and $\overline z$. Note that $1/z$ is
 on the same circle. Denote by $\lambda_z$ the arc on $C_z$ with endpoints
 $\overline z$ and $1/z$ which does not contain $t$ and $1/t$. We shall show that
 $$
 M_z\hbox{\ consists of \ } \lambda_z \hbox{\ and of the segment joining\ }\overline z\hbox{\ and\ } 1/z.
 \eqno (8.1)
 $$
 We have $L\cap(\{ z\} \times\C )=\{ (z,R^2/z)\colon\ |z|\leq R\leq1\}$ which
 is the segment joining $(z,\overline z)$ and $(z,1/z)$. To
 find what $N\cap (\{ z\}\times\overline \C\}$ is, we recall first that
 $$\eqalign{
    \Lambda(T,\sqrt\Tt) = \bigl\{ (z,w)\colon\ w=T+&{{\Tt}\over{z-T}},\cr
                                                 &|z-T|<\sqrt\Tt\bigr\}.\cr}
    $$
    So we must determine $\{ w(T)\colon\ 0\leq T<T(z)\}$ where
    $$
    w(T) = T+{{\Tt}\over{z-T}} = {{Tz-tT-T/t+1}\over{z-T}}
    $$
    and where $T(z)$ is such that $w(T(z))=\overline z$. The set $\{ w(T)\colon\ T\in\R\}$
    is clearly a circle. We have $w(t)=t,\  w(1/t)=1/t$ and
 $w(0)=1/z$. So this circle is $C_z$ described above and it follows that when
 $0\leq T\leq T(z)$ then $y(T)$ is on the arc $\lambda_z$. This proves (8.1).

 We now look at $M_x$ where $x\in \cS$ is real, By (7.4) we know that if $(x,y)\in\cT_R$ then, as
 $R$ increases from $x$ to $1$,\ $y$\ increases from
 $x$ to $1/x$. By (7.6) we know that if $(x,y)\in\cS_T\cup\{ (T,\infty)\} $ then,
 as T increases from $0$ to $x,\ \ y$ increases from $1/x$ to $\infty$ and as
 $T$ increases from $x$ to $T_0$,\ \ $y$ increases from $-\infty$ to $x$. This shows that
 $$
 \hbox{if\ } x\in\cS\hbox{\ is real then\ }M_x\hbox{\ is the real axis }.
 \eqno (8.2)
 $$

 For each $z\in\cS$, let $D_z$ be the domain bounded by $M_z$, The domains
 $D_z$ change continuously with $z\in\cS\setminus\R$ and as
 $z$ approaches a point $a\in\bD\setminus\R$ they shrink to the point $\overline a$.

 At this point we are precisely in the situation described in Section 4 of [G6].
 Repeating word by word the part of the proof there
 we use an argument of H.Lewy [L] as generalized by H.Rossi [R], to prove that
 for each $z\in\cS\setminus \R$ the function $w\in\Phi (z,w)$, defined on $M_z$,
 extends holomorphically through $D_z$. 
 
 Recall that
 $\{ (\lambda,\infty)\colon\ 0<\lambda<\eta\}\subset M$ and that $\Phi$ is continuous
 on $M$. Given $\tau ,\ 0<\tau<\eta$, we show that $\Phi $ is
 constant on $\{\tau\}\times M_\tau$. To see this, recall that $\eta<t$ and fix $\tau,\ 0<\tau<\eta$.  
 Observe that for small $\omega >0$, \ $M_{\tau+i\omega}$
 are simple closed curves bounding $D_{\tau+i\omega}$ which depend continuously
 on $\omega $ and, as domains in $\overline\C$ converge to the
 halfplane $\Im\z<0$ as $\omega $ tends to $0$. Since for each small $\omega $
 the function $\z\mapsto\Phi (\tau+i\omega,\z )$ extends from $M_{\tau + i\omega}$
 holomorphically through $D_{\tau+i\omega}$,
 the continuity of $\Phi $ on $M$ implies that $s\mapsto\Phi (\tau, s)$ has a
 bounded continuous extension from $\R$ to the halfplane $\Im\z\leq 0$
 which is holomorphic  on $\Im\z <0$. Repeating the reasoning with $\omega<0$
 we see that $s\mapsto \Phi (\tau, s)$ has
 a bounded continuous extension from $\R$ to the halfplane $\Im\z\geq 0$
 which is holomorphic  on $\Im\z >0$. Thus, $s\mapsto\Phi (\tau, s)$ has a
 bounded holomorphic extension to $\C$, which, by the Liouville theorem, must be constant.
 Thus, for each $\tau,\ 0<\tau<\eta$, the holomorphic 
 extensions of $\varphi $ from all cirles in our families that surround $\tau $, coincide
 at $\tau $. This implies that
 $\varphi $ is holomorphic in a neighbourhood of the segment $(0,\eta )$ and it is then easy to
 see that the analyticity propagates along circles so $\varphi $ is holomorphic on $\D$.
 This completes the proof of Lemma 5.1. and also proves Theorem 1.3. The proof of Lemma 3.2 is thus complete.
 This completes the proof of Theorem 1.2.
 \vskip 4mm
 \bf 9.\ Higher dimensions\rm
 \vskip 2mm
 After the results presented above were obtained Mark Agranovsky told the author that
 he has proved Theorem 1.1. for real analytic functions and its generalization to
 higher dimensions and then put his results on the web [A1, A2]. In a later version of
 [A1] he showed that it is very easy to generalize his result to higher dimensions and that,
 surprisingly, only complex lines through two points in $\overline B$ suffice in any dimension.
 We want to repeat the elegant simple argument that he used and apply it in our case.
 \vskip 2mm
 \noindent\bf COROLLARY 9.1 \ \it Let $B$ be the open unit ball in $\C^N,\ N\geq 2$, and
 assume that $a, b\in\ \C^N,\ a\not=b$, are
 such that the complex line containing
 $a$ and $b$ meets $B$ and such that $<a|b>\not=1$.
 Suppose that a function $f\in\C^\infty (bB)$ extends holomorphically along
 every complex line in $\L (a)\cup\L (b)$. Then $f$ extends holomorphically through $B$.\rm
 \vskip 1mm 
 \bf Proof.\ \rm Using Moebius transforms  we may, with no loss of generality assume that
 $\Lambda (a,b)$, the complex line through $a$ and $b$, contains the origin. Choose $c\in\Lambda(a,b)\cap B,\ c\not= 0$.
 By Theorem 1.2 the
 function $f|(\Sigma\cap bB)$ extends holomorphically through $\Sigma\cap B$ for every complex two-plane
 $\Sigma$ containing $\Lambda (a,b)$. So, if a complex line $L\in\L (0)$ meets a complex line $E\in\L (c)$
 then both $L$ and $E$ are contained in such a two-plane
 which implies that $f$ extends holomorphically along both $E$ and $L$ and that the extensions are the same
 at $E\cap L$. Thus, all such
 holomorphic extensions along complex lines in
 $\L (0)\cup\L (c)$ define an extension $\tilde f$ of $f$ to $\overline B$ which is holomorphic
 on $L\cap B$ for each $L\in\L (0)$ and holomorphic on $E\cap B$ for each $E\in\L (c)$. Expressing
 $\tilde f| E\cap B$ as the Cauchy integral of $f|E\cap bB$ for each $E\in\L (c)$,
 the fact that $f\in C^\infty (bB)$ implies that $\tilde f$ is of class $C^\infty $ in a neighbourhood
 of the origin which, by a theorem of F.\ Forelli [Ru]
 implies that $\tilde f$ is holomorphic on $B$, which completes the proof.
 \vskip 1mm
 Very recently L.\ Baracco [B] showed that for real analytic functions on $bB$ only one point in the boundary
 suffices, that is,
 if $a\in bB$ then $\L (a)$ is a test family for holomorphic extendibility for real analytic functions on $bB$.
 \vskip 4mm
 \bf 10.\ An open question
 \vskip 2mm
 \rm 
 Let us conclude by formulating our initial question for general domains:
 \vskip 1mm
 \noindent\bf QUESTION 10.1\ \it Let $a,b\in\C^2,\ a\not=b$, and let $D\subset \C^2$ be a bounded domain with
 smooth boundary. Assume that a continuous function $f$ on $bD$
 extends holomorphically along every complex line $L\in\L (a)\cup\L (b)$ that meets $D$. Must $f$ extend
 holomorphically through $D$? \rm
 \vskip 1mm
 \noindent We have seen that if $D$ is a ball the answer is no even in the case when $a,b\in D$ and to get
 a positive answer one must assume that $f$ is infinitely smooth.
 The author believes that one must require smoothness only in very special cases:
 \vskip 1mm
 \noindent\bf CONJECTURE\ \it For a generic domain $D$ the
 answer to Question 10.1 is positive. \rm
 \vskip 1mm
 \noindent M.\ Lawrence has proved a result of this kind for arbitrary small perturbations of the ball [La].
 \vskip 4mm
 \noindent\bf ACKNOWLEDGEMENTS\ \rm The author is grateful to Yum-Tong Siu whose lecture at a conference in 
 Seattle in 1998
 raised Question 10.1 and made the author conjecture that Question 1.1 has a positive answer for functions of
 class $C^\infty$. The author is grateful to Mark Agranovsky for stimulating discussions. A part of this paper was
 written in November 2009 during the author's stay at the
 E.\ Schroedinger Institute for Mathematical Physics in Vienna. He wishes to thank Fritz Haslinger, Bernard Lamel
 and Emil Straube for the invitation.
 \vskip 3mm
 This paper was supported in part by the ministry of Higher Education, Science and Technology of Slovenia through the research
 program Analysis and Geometry, Contract No. P1-02091 (B).
\vskip 15mm
\centerline{\bf REFERENCES}
\vskip 5mm
\noindent [A1]\ M.\ L.\ Agranovsky:\ Holomorphic extension from the unit
sphere in $C^n $ into complex lines passing through a finite set.

\noindent http://arxiv.org/abs/0910.3592
\vskip 2mm
\noindent [A2]\ M.\ L.\ Agranovsky:\ Characterization of polyanalytic functions by meromorphic extensions into
chains of circles.

\noindent http://arxiv.org/abs/0910.3578
\vskip 2mm
\noindent [AG]\ M.\ L.\ Agranovsky and J.\ Globevnik:\ Analyticity
on circles for rational and real-analytic functions of two real variables.

\noindent J.\ d'Analyse Math.\ 91 (2003) 31-65
\vskip 2mm
\noindent [B]\ L.\ Baracco: Holomorphic extension from the sphere to the ball.

\noindent http://arxiv.org/abs/0911.2560
\vskip 2mm
\noindent [G1]\ J.\ Globevnik;\ On holomorphic extensions from spheres in $C^2$

\noindent Proc.\ Roy.\ Soc.\ Edinb.\ 94A (1983) 113-120
\vskip 2mm
\noindent [G2[\ J.\ Globevnik:\ A family of lines for testing holomorphy in the ball of $C^2$.

\noindent Indiana Univ. Math. J.\ 36 (1987) 639-644
\vskip 2mm
\noindent [G3] \ J.\ Globevnik:\ Holomorphic extensions from open families of circles.

\noindent Trans.\ Amer.\ Math.\ Soc.\ 355 (2003) 1921-1931
\vskip 2mm
\noindent [G4]\ J.\ Globevnik:\ Analyticity on families of circles.

\noindent Israel J.\ Math. 142 (2004) 29-45
\vskip 2mm
\noindent [G5] \ J.\ Globevnik:\ Analyticity on translates of a Jordan curve.

\noindent Trans.\ Amer.\ Math.\ Soc. 359 (2007) 5555-5565
\vskip 2mm
\noindent [G6]\ J.\ Globevnik:\ Analyticity of functions analytic on circles

\noindent Journ.\ Math.\ Anal.\ Appl. 360 (2009) 363-368
\vskip 2mm
\noindent [KM]\ A.\ M.\ Kytmanov, S.\ G.\ Myslivets:\ On families of complex lines sufficient for holomorphic extension.

\noindent Mathematical Notes 83 (2008) 500-505 (translated from Mat.Zametki 83 (2008) 545-551
\vskip 2mm
\noindent [La]\ M.\ G.\ Lawrence:  Hartogs' separate analyticity theorem for CR functions.

\noindent Internat.\ J.\ Math.\ 18 (2007) 219--229.
\vskip 2mm
\noindent [L]\ H.\ Lewy: On the local character of the
solutions of an atypical linear differential equation in three
variables and a related theorem for regular functions of two complex variables.

\noindent Ann. of Math.\ 64 (1956) 514-522
\vskip 2mm
\noindent [R]\ H.\ Rossi: A generalization of a theorem of Hans Lewy.

\noindent Proc.\ Amer.\ Math.\ Soc.\ 19 (1968) 436-440
\vskip 2mm
\noindent [Ru]\ W.\ Rudin:\ \it Function theory in the unit ball of $C^n$.\rm

\noindent Springer, Berlin-Heidelberg-New York 1980
\vskip 2mm
\noindent [T1]\ A.\ Tumanov: A Morera type theorem in the strip.

\noindent Math.\ Res.\ Lett.\ 11 (2004) 23-29
\vskip 2mm
\noindent [T2]\ A.\ Tumanov: Testing analyticity on circles.

\noindent Amer.\ J.\ Math.\ 129 (2007) 785-790

\vskip 8mm
\noindent Institute of Mathematics, Physics and Mechanics

\noindent University of Ljubljana, Ljubljana, Slovenia

\noindent josip.globevnik@fmf.uni-lj.si

\bye